# Pranab Kumar Sen: Life and works

**N. Balakrishnan[1], Edsel A. Peña[2] and Mervyn J. Silvapulle** [3]

*McMaster University, University of South Carolina and Monash University*

**Abstract:** In this article, we describe briefly the highlights and various accomplishments in the personal as well as the academic life of Professor Pranab Kumar Sen.

Pranab Kumar Sen (born: November 7, 1937) had his school and college education (B.Sc (1955), M.Sc. (1957) and Ph.D. (1962), all in Statistics) from Calcutta University, where he served as a lecturer in the Post-graduate Department of Statistics (1961–64) before proceeding to the University of California, Berkeley (1964–65) as a visiting faculty member in Statistics. He moved to the University of North Carolina (UNC) at Chapel Hill in the Fall of 1965, and since then, he has been there (visiting 1965–67, Associate Professor 1967–70 and Professor 1970–). He is currently Cary C. Boshamer Professor of Biostatistics (1982–), and Professor of Statistics and Operations Research (1988–) at UNC. He was Richard Metron Guest Professor at the Albert-Ludwig University, Freiburg, Germany (1974–75) and also the Eugene Lukacs Distinguished Professor, Bowling Green State University, Ohio (1997 Spring). In 1993 he was made an Adjunct (life-long) Professor, Indian Statistical Institute. He has travelled extensively and visited many Universities all over the world, and quite frequently to University of Sao Paulo, Brazil, Institute of Statistical Science, Academia Sinica ( and National Health Research Institute), Taipei, Taiwan, Charles University, Prague, Czech Republic, and the Indian Statistical Institute.

Professor Sen was elected Fellow of the Institute of Mathematical Statistics in 1968, Fellow of the American Statistical Association in 1969, and Elected Member of the International Statistical Institute in 1973. He was a NSF-CBMS Lecturer (1983) at the University of Iowa, S.N. Roy Memorial Lecturer, Calcutta University, 1976, Fifth Special Invited Lecturer at the Institute of Statistical Science, Academia Sinica, Taipei (2001), and Platinum Jubilee Lecturer, Indian Statistical Institute (2007). He was awarded the Boltzman Medal, Czech Union of Mathematicians and Physicists (CUMP) in 1988, and in 1998, the CUMP Commemoration Medal for outstanding contributions to Statistics and Probability Theory. At UNC, in 1996, he was awarded the McGavran Teaching Excellence Award. In 2002, the American Statistical Association presented the Senior Noether Award to him for life-long achievements in research and teaching (nonparametrics).

Professor Sen has supervised more than 80 doctoral students (1969–2007) in Statistics, Biostatistics, and Psychometry, all at UNC. He has (co-)authored and

[1]Department of Mathematics and Statistics, McMaster University, Hamilton, Ontario, Canada L8S 4K1, e-mail: bala@univmail.cis.mcmaster.ca

[2]Department of Statistics, University of South Carolina, Columbia, SC 29208, USA, e-mail: pena@stat.sc.edu

[3]Department of Econometrics and Business Statistics, Monash University, Caulfield East, Australia 3145, e-mail: mervyn.silvapulle@buseco.monash.edu.au

(co-)edited 22 books and monographs in Statistics, Biostatistics, Stochastic Processes, and related fields. Besides, he has more than 600 publications in Statistics, Probability Theory, Stochastic Processes, and Biostatistics in leading journals in these areas. He has served on the Editorial Boards of Journal of Multivariate Analysis (1972–78), Journal of the American Statistical Association (1973–78), Communications in Statistics, Theory and Methods (1972–82), Journal of Statistical Planning and Inference (1979–90, 1995–2006), Brazilian Journal of Probability and Statistics (1988–), Journal of Nonparametric Statistics (1991–94), Lifetime Data Analysis (2005–), Metron (2005–07), and Scientiae Mathematicae Japonicae (2003–). He was the Founding (joint) Editor of Sequential Analysis (1982–95), also of Statistics and Decisions (1982–2002). He is the Chief Editor, Sankhya (2007–2009).

Sen's research impact has mostly been in the original and innovative developments of nonparametrics and asymptotics in a broad spectrum of Statistical (and Biostatistical) Science, encompassing both theory and applications in interdisciplinary fields. In his doctoral work, he developed the jackknife theory for U-statistics (1960), and with S.K. Chatterjee (1964), laid down the foundation of multivariate nonparametrics; the field was evolutionary in the late 1960's and culminated with the monograph of Puri and Sen (1971). In 1970–74, Malay Ghosh and Sen did most innovative research in sequential nonparametrics, and with S. K. Chatterjee (1972–74), he developed the theory of time-sequential nonparametrics, both culminating in his 1981 monograph and 1983 NSF-CBMS Lectures. With Madan Puri, he had a second monograph on nonparametrics for general linear models (1985). At that time, he started working on robustness with Jana Jurečková and their long collaboration led to their 1996 monograph on this field of research. Sequential analysis and applications in clinical trials captured much of his interest during the 1980's and 1990's, and the 1997 monograph relates to some of these developments. Pitman Closeness measure has been a major point of interest for his research, and the Keating–Mason–Sen (1993) monograph reveals a lot of their work along with others. At the intermediate level, the Sen–Singer (1993) monograph served as a good text/reference for many graduate students in Biostatistics and Applied Satatistics fields. Two other notable monographs with late Jaroslav Hajek and Z. Šidák (1999) and Mervyn Silvapulle (2005) have received considerable attention from Statistics Communities all over the world. During the past 15 years, Sen is more inclined to research in Biostochastics, Environmetrics, and Bioinformatics and continuing basic methodological work to bridge the gap between theory and applications. Beyond Parametrics is a phrase that Sen coined to encompass the entire field of Nonparametrics and Semiparametrics (including Bayesian methods) that captures better the reality of interdisciplinary research and serves as the guiding post for statistical innovations.

We have taken the liberty of including the following verse Sen composed a couple of years ago.

## My chancy life as a statistician

Pranab Kumar Sen
University of North Carolina, Chapel Hill

> Is there any human life without the chance in manifestation?
> Could it be any less chancy for a converted (bio-)statistician?
> When I was born, the chances were that I might not even survive,

not to speak of the fact that I would be in a stochastic blow-pipe!
Although not from an affluent family, the first ten years swept smoothly.
Then the first throw of a dice was the demise of my father, abruptly.
Crippled with the political partition of Bengal (albeit independence),
the chance for survival became too much relatives' dependence.
The last four years of the high school was without any glamour,
my mind was set for more attractions on the roadside games-parlours.
My poor mother gave up all the hopes for any smooth matriculation,
yet by a second throw of a dice, I topped the list of school graduation.
I was trying to optimize my chances of getting into a medical school,
alas, being underaged, I was denied chancy entrance to the physicians' pools.
In disgust, I ran to the (Calcutta) Presidency College courtyards,
and grabbed a spot in both mathematics and physics honours charts.
A third throw of a dice: came an old friend in utter surprise,
and asked me to switch to statistics major on his chancy advise.
The family financial burden was about to knock me out of college,
but the chancy success in the B.Sc. examination, opened the PG page.
With the march of countless dice, I could step outside the masters,
and embarked on a chancy ambition of wearing doctorate feathers.
Flanked by my chancy success, I stood at the ISICAL doorsteps,
in quest of an assistantship, for which I was even aptitude tested!
Another throw of a dice: I was not selected. Not because of merit,
nor of underage, but for not enough money to cover my family need.
I turned my face back to my alma mater, and locked my chance
of being in predoctoral training program in thrilled enhance.
Those six years in the Calcutta University were the paradise of my life,
I could have a glance of the statistics interface in broad delight.
I got a chance to visit the University of California, Berkeley, in 1962,
however, family situation kept me anchoring to Calcutta for a year or two.
A shower of chance in 1964: via Berkeley, we came to Chapel Hill,
a small town, home for the rest of my life, in a chancy exile.
I never dreamed of becoming a statistician, not to speak of an applied one,
yet the nonparametric I specialized, chanced me to this open direction.
Back in 1961 when I was asked to introduce a course in bioassay,
the odds were in my favour for pouring nonparametrics in this say.
There, Shoutir Chatterjee and I embarked on multivatiate nonparametrics,
a chancy step that extended all the way to Chapel Hill in all peroratives.
Fortune to have Madan Puri on my side in incessant youthful delight,
to carry the rank-permutation banner; a combination eventually out of sight.
Time-sequential as well as sequential nonparametrics evolved soon after,
Shoutir Chatterjee and Malay Ghosh set the chancy collaboration to ponder.
At that time, I was more rank (analysis) minded, albeit with some robustness,
and had a chance to set collaboration with Jana Jurečková in a long process.
Got a chance to traverse preliminary test inference to a shrunken base,
with the company of Saleh, Kubokawa and Keating for the Pitman closeness.
In quest of this, the finite sample approaches led to the asymptotics,
and Julio Singer joined me to build a sand-castle in simple acquities.
Could not refuse an appealing call from Mrs. Hájek and Zbyněk Šidák
to restore and update a classical text on the theory of test based on rank;
it took a while to add new material in a unifying but still original style,
all for the respect to late Jaroslav Hájek, not for any personal stake.
Life is full of constraints, contrary to our expectation,
in statistical inference too, constraints are no exception.
I was delighted to have Mervyn Silvapulle in active cooperation,
taking the lead in a constraint inference monograph in resolution.

The agony of being denied medical college entrance chanced me a signal:
stick to biostatistics to fathom out chance mysteries of a clinical trial.
Soon after bioinformatics and genomics rocked the field with data mining,
and I have the chance to come up with biostochastics in the offing.
Let me cherish the chance beauty of stochastics in applicable methodology,
and gauge knowledge discovery and data mining, in the chance ideology.

## RESEARCH MONOGRAPHS AND ADVANCED TEXT BOOKS

(1) *Nonparmetric Methods in Multivariate Analysis* (co-author: M.L. Puri), John Wiley and Sons, New York, 197l.

(2) *Sequential Nonparametrics: Invariance Principles and Statistical Inference.* John Wiley and Sons, New York, 1981.

(3) *Contributions to Statistics: Essays in Honor of Norman L. Johnson* (edited). North-Holland, Amsterdam, 1983.

(4) *Handbook of Statistics, Volume 4: Nonparametric Methods* (co-editor: P.R Krishnaiah), North-Holland, Amsterdam, 1984.

(5) *Biostatistics: Statistics in Biomedical, Public Health and Environmental Sciences* (edited). North-Holland, Amsterdam, 1985.

(6) *Nonparametric Methods in General Linear Models* (co-author: M.L. Puri), John Wiley and Sons, New York, 1985.

(7) *Theory and Applications of Sequential Nonparametrics.* SIAM Publication (CBMS NSF), Philadelphia, 1985.

(8) *Goodness of Fit* (co-editor with Pal Revesz and Karoly Sarkadi), North-Holland, Amsterdam, l986.

(9) *Handbook of Sequential Analysis* (co editor: B.K. Ghosh), Marcel Dekker, N.Y. 1991.

(10) *Pitman's Measure of Closeness: A Comparison of Statistical Estimators* (co-authors: J.P. Keating and R.L. Mason), SIAM, Philadelphia, 1993.

(11) *Large Sample Methods in Statistics: An Introduction with Applications* (co-author: J.M. Singer), Chapman and Hall, 1993.

(12) *Order Statistics and Nonparametrics: Theory and Applications* (co-editor: I. A. Salama), North-Holland, Amsterdam, 1992.

(13) *Stochastic Processes, a Festschrift in honour of Gopinath Kallianpur* (co-editors: S. Cambanis, J.K. Ghosh and R.L. Karandikar), Springer-Verlag, New York, 1993.

(14) *Collected Works of Wassily Hoeffding* (co-editor: N. I. Fisher), Springer-Verlag, 1993.

(15) *Statistical Theory and Applications: Papers in Honor of Herbert A. David* (co-editors H.N. Nagaraja and D. Morrison), Springer-Verlag, New York, 1995.

(16) *Robust Statistical Procedures: Asymptotics and Inter-relations* (co-author: J. Jurečková), John Wiley and Sons, New York, 1996.

(17) *Sequential Estimation* (co-authors: M. Ghosh and N. Mukhopadhyay), John Wiley and Sons, New York, 1997.

(18) *Theory of Rank Tests, 2nd Edition* (co-authors: J. Hájek and Z. Šidák). Academic Press, United Kingdom, 1999.

(19) *Handbook of Statistics, Volume 18: Bioenvironmental and Public Health Statistics,* (co-editor: C. R. Rao). Elsevier, Amsterdam, 2000.

(20) *Perspectives in Statistical Sciences,* (co-editors: A. Basu, J. K. Ghosh, P. K. Sen and B. K. Sinha. Oxford India, Delhi, 2000.

(21) *Constrained Statistical Inference: Inequality, Order, and Shape Restrictions,* (co-author : M. J. Silvapulle), John Wiley and Sons, New York, 2005.

(22) *Excursions in Biostochastics: Biometry to Biostatistics to Bioinformatics,* Invited Lecture Notes, No. 5, Academia Sinica, Taipei, Taiwan, 2004.

## PH.D. GUIDANCE

Served as the adviser and supervised the doctoral dissertations of the following persons (all at the University of North Carolina, Chapel Hill):

| | | | |
|---|---|---|---|
| 1. R.R. Obenchain | (Mathematical Statistics) | 1969 |
| 2. M. Ghosh | (Mathematical Statistics) | 1969 |
| 3. T.M. Gerig | (Mathematical Statistics) | 1971 |
| 4. K. Dutta | Mathematical Statistics) | 1971 |
| 5. G.W. Williams | (Biostatistics) | 1972 |
| 6. R.F. Woolson | (Biostatistics) | 1973 |
| 7. M.R. Mahmoud | (Biostatistics) | 1973 |
| 8. Y. Hochberg | (Biostatistics) | 1974 |
| 9. K.L. Monti | (Biostatistics) 1975 | |
| | (co-adviser: D. Quade) | |
| 10. H. Majumdar | (Biostatistics) | 1976 |
| 11. C. Silva | (Biostatistics) | 1977 |
| | (co-adviser: D. Quade) | |
| 12. J.C. Gardiner | (Mathematical Statistics) | 1978 |
| | (co-adviser: G. Simons) | |
| 13. Y. Tsong | (Mathematical Statistics) | 1979 |
| 14. V.M. Chinchilli | (Mathematical Statistics) | 1979 |
| 15. F.E. Harrell, Jr. | (Biostatistics) | 1979 |
| 16. A.N. Sinha | (Biostatistics) | 1979 |
| 17. E.R. DeLong | (Biostatistics) | 1979 |
| 18. S.I. Bangdiwala | (Biostatistics) | 1980 |
| 19. Y.C. So | (Mathematical Statistics) | 1981 |
| 20. I.N. Agung | (Biostatistics) | 1981 |
| 21. K. Kouri | (Biostatistics) | 1981 |
| 22. G. Feeney | (Biostatistics) | 1982 |
| | (co-adviser: M. Symons) | |
| 23. M.N. Boyd | (Biostatistics) | 1982 |
| 24. J. Singer | (Biostatistics) | 1983 |

| 25. D.L. Hawkins | (Biostatistics) | 1983 |
|---|---|---|
| 26. R.A. Smith | (Mathematical Statistics) | 1983 |
| | (co-adviser: W. Hoeffding) | |
| 27. Y.C. Yuan | (Biostatistics) | 1984 |
| 28. R.W. Falk | (Biostatistics) | 1985 |
| 29. D. Hoberman | (Biostatistics) | 1986 |
| 30. A.R. Karmous | (Biostatistics) | 1986 |
| 31. M.T. Tsai | (Mathematical Statistics) | l987 |
| 32. G.R. Jerdack | (Biostatistics) | l987 |
| 33. D.C. Trost | (Biostatistics) | 1988 |
| 34. D. Sengupta | (Mathematical Statistics) | 1988 |
| 35. Rick Williams | (Biostatistics) | 1988 |
| 36. C.Y. Wada | (Biostatistics) | 1988 |
| 37. S.A. Murphy | (Mathematical Statistics) | 1989 |
| 38. J.I. Crowell | (Mathematical Statistics) | 1990 |
| 39. J.J. Ren | (Mathematical Statistics) | 1990 |
| 40. L.J. Edwards | (Biostatistics) | 1990 |
| 41. Nelson Oliveira | (Biostataistics) | 1992 |
| | (co-advisor: D. Quade) | |
| 42. Subha Das | (Mathematical Statistics) | 1993 |
| 43. Ming Zhang | (Mathematical Statistics) | 1993 |
| 44. Ralph DeMasi | (Biostatistics) | 1994 |
| | (co-advisor: B. Qaqish) | |
| 45. S. Munoz | (Biostatistics) | 1995 |
| | (co-advisor: S. I. Bangdiwala) | |
| 46. A. Pedroso de Lima | (Biostatistics) | 1995 |
| 47. Habib El-Moalem | (Biostatistics) | 1995 |
| 48. Moh. A. Chaudhary | (Biostatistics) | 1995 |
| 49. Kentaro Hayashi | (Psychometry) | 1996 |
| 50. Zhenwei Zhou | (Mathematical Statistics) | 1996 |
| 51. Pai-Lien Chen | (Biostatistics) | 1996 |
| 52. Solange Andreoni | (Biostatistics) | 1996 |
| 53. Stoffel Moeng | (Biostatistics) | 1996 |
| | (co-advisor: C.M. Suchindran) | |
| 54. Ho Kim | (Biostatistics) | 1996 |
| 55. Lin X. Clegg | (Biostatistics) | 1997 |
| | (co-advisor: Jianwen Cai) | |
| 56. Mario Chen Mok | (Biostatistics) | 1997 |
| 57. Maha Karnoub | (Biostatistics) 1997 | |
| | (co-advisor: Francoise Seillier-Moiseiwitsch) | |
| 58. Hildete P. Pinheiro | (Biostatistics) | 1997 |
| | (co-advisor: Francoise Seilier-Moiseiwitsch) | |
| 59. H.-C. Tien | (Biostatistics) | 1998 |
| 60. Michael Marion | (Mathematical Statistics) | 1998 |
| 61. J. E. MacDougal | (Biostatistics) | 1999 |
| 62. Ming Zhong | (Biostatistics) | 1999 |
| | (co-advisor: J. Cai) | |
| 63. Antonio J. Sanhueza | (Biostatistics) | 2000 |
| 64. H. Chakravarty | (Dr.P.H. Biostatistics) | 2000 |
| | (co-advisor: R. Helms) | |
| 65. Karen Kessler | (Biostatistics) | 2000 |
| 66. Dubois Bowman | (Biostatistics) | 2000 |
| | (co-advisor: P. Stewart) | |
| 67. Joe Galanko | (Biostatistics) | 2000 |
| | (co-advisor: D. Quade) | |

| 68. Chula Komoltri | (Dr. P.H. Biostatistics) | 2001 |
| | (co-advisor: K. Bangdiwalla) | |
| 69. Jianmin Wang | (Biostatistics) | 2001 |
| | (co-advisor: C. M. Suchindran) | |
| 70. Kouros Owzar | (Mathematical Statistics) | 2002 |
| 71. Stephanie Cano | (Biostatistics) | 2003 |
| | (co-advisor: C. M. Suchindran) | |
| 72. Tejas Desai | (Biostatistics) | 2003 |
| 73. Lan Kong | (Biostatistics) | 2003 |
| | (co-advisor: Jianwen Cai) | |
| 74. Szu-Yun Leu | (Biostaistics) | 2003 |
| 75. Inkyung Jung | (Biostatistics) | 2004 |
| 76. Lily Wang | (Biostatistics) | 2004 |
| 77. George Capuano | (Biostatistics) | 2005 |
| 78. Munni Begum | (Dr.P.H., Biostatistics) | 2005 |
| 79. George Luta | (Biostatistics) | 2006 |
| | (co-advisor: Gary Koch) | |
| 80. Munsu Kang | (Biostatistics) | 2007 |

## Selected publications of Pranab K. Sen